\documentclass[10pt,twoside,a4paper]{scrartcl}

\usepackage[left=3cm,right=4cm,top=3cm,bottom=4cm,includeheadfoot]{geometry}
\usepackage[english]{babel}
\usepackage[utf8]{inputenc}
\usepackage{fancyhdr}
\usepackage{amsmath, amssymb, latexsym,
   epsfig, rotating, fancyhdr, amsthm, pifont}
\usepackage{amsmath}
\usepackage{verbatim} 

\newcommand{\foot}{\footnote}

\newcommand{\ld}{\ensuremath{,\ldots,}}
\newcommand{\ssq}{\ensuremath{\subseteq}}
\newcommand{\smin}{\ensuremath{\setminus}}
\newcommand{\eps}{\ensuremath{\varepsilon}}

\newcommand{\htop}{\ensuremath{h_{\mathrm{top}}}}

\newcommand{\ind}{\ensuremath{\mathbf{1}}}

\newcommand{\diam}{\ensuremath{\mathrm{diam}}}

\newcommand{\kreis}{\ensuremath{\mathbb{T}^{1}}}

\newcommand{\torus}{\ensuremath{\mathbb{T}^2}}

\newcommand{\nfolge}[1]{\ensuremath{(#1)_{n\in\mathbb{N}}}}

\newcommand{\twovector}[2]{\ensuremath{\left(\begin{array}{c} #1 \\
        #2 \end{array}\right)}}

\newcommand{\alphlist}{\begin{list}{(\alph{enumi})}{\usecounter{enumi}\setlength{\parsep}{2pt}
      \setlength{\itemsep}{1pt} \setlength{\topsep}{5pt}
      \setlength{\partopsep}{3pt}}}
\newcommand{\arablist}{\begin{list}{(\arabic{enumi})}{\usecounter{enumi}\setlength{\parsep}{2pt}
          \setlength{\itemsep}{1pt} \setlength{\topsep}{5pt}
          \setlength{\partopsep}{3pt}}}
\newcommand{\romanlist}{\begin{list}{(\roman{enumi})}{\usecounter{enumi}\setlength{\parsep}{2pt}
              \setlength{\itemsep}{1pt} \setlength{\topsep}{5pt}
              \setlength{\partopsep}{3pt}}}
\newcommand{\Romanlist}{\begin{list}{(\Roman{enumi})}{\usecounter{enumi}\setlength{\parsep}{2pt}
              \setlength{\itemsep}{1pt} \setlength{\topsep}{5pt}
              \setlength{\partopsep}{3pt}}}
\newcommand{\bulletlist}{\begin{list}{$\bullet$}{\setlength{\parsep}{2pt}
                \setlength{\itemsep}{1pt} \setlength{\topsep}{5pt}
                \setlength{\partopsep}{3pt}\setlength{\leftmargin}{15pt}}} 
\newcommand{\Alphlist}{\begin{list}{(\Alph{enumi})}{\usecounter{enumi}\setlength{\parsep}{2pt}
      \setlength{\itemsep}{1pt} \setlength{\topsep}{5pt}
      \setlength{\partopsep}{3pt}}}
 \newcommand{\listend}{\end{list}}

\newcommand{\N}{\ensuremath{\mathbb{N}}} 
\newcommand{\R}{\ensuremath{\mathbb{R}}}
\newcommand{\Z}{\ensuremath{\mathbb{Z}}}
\newcommand{\Q}{\ensuremath{\mathbb{Q}}}

\newcommand{\cA}{\mathcal{A}}

\newcommand{\cC}{\mathcal{C}}
\newcommand{\cD}{\mathcal{D}}

\newcommand{\cT}{\mathcal{T}}

\newcommand{\nLim}{\ensuremath{\lim_{n\rightarrow\infty}}}

\newcommand{\kLim}{\ensuremath{\lim_{k\rightarrow\infty}}}

\newcommand{\inergsum}{\ensuremath{\sum_{i=0}^{n-1}}}

\newcommand{\ntel}{\ensuremath{\frac{1}{n}}}

\newcommand{\halb}{\ensuremath{\frac{1}{2}}}

\newcommand{\viertel}{\ensuremath{\frac{1}{4}}}

\newtheorem{theorem}{Theorem}[section]
\newtheorem{corollary}[theorem]{Corollary}

\newtheorem{lemma}[theorem]{Lemma}

\theoremstyle{definition}

\newtheorem*{example*}{Example}
\newtheorem{remark}[theorem]{Remark}

\newcommand{\fsc}{\mathrm{{ac}}}

\newcommand{\Per}{\ensuremath{\mathrm{Per}}}

\newcommand{\hpow}{\ensuremath{h_{\mathrm{pow}}}}
\newcommand{\hmod}{\ensuremath{h_{\mathrm{mod}}}}

\newcommand{\hpowup}{\ensuremath{\overline{h}_{\mathrm{pow}}}}
\newcommand{\hpowlow}{\ensuremath{\underline{h}_{\mathrm{pow}}}}

\newcommand{\hmodup}{\ensuremath{\overline{h}_{\mathrm{mod}}}}
\newcommand{\hmodlow}{\ensuremath{\underline{h}_{\mathrm{mod}}}}

\title{\Large\textsc{Some remarks on modified power entropy}}
\author{\normalsize M. Gr\"oger \thanks{Department of Mathematics,
    Universit\"at Bremen, Germany. Email: {\tt
      groeger@math.uni-bremen.de}} \and \normalsize T. J\"ager
  \thanks{Department of Mathematics, TU Dresden, Germany. Email: {\tt
      Tobias.Oertel-Jaeger@tu-dresden.de}} }

\begin{document}
\renewcommand\dagger{**}
\maketitle

\begin{abstract}\small The aim
  of this note is to point out some observations concerning modified power
  entropy of $\Z$- and $\N$-actions. First, we provide an elementary example
  showing that this quantity is sensitive to transient dynamics, and therefore
  does not satisfy a variational principle. Further, we show that modified power
  entropy is not suitable to detect the break of equicontinuity which takes
  place during the transition from almost periodic to almost automorphic minimal
  systems. In this respect, it differs from power entropy and amorphic
  complexity, which are two further topological invariants for zero entropy
  systems (`slow entropies'). Finally, we construct an example of an irregular
  Toeplitz flow with zero modified power entropy.
\end{abstract}

\section{Introduction}
\label{Intro}
Given a continuous map $f:X\to X$ on some compact metric space
$(X,d)$, the {\em Bowen-Dinaburg metrics} are given by 
\begin{equation*}
  d_n^f(x,y) \ = \ \max_{i=0}^{n-1} d(f^i(x),f^i(y)) \ .
\end{equation*}
If $S_n(f,\delta)$ denotes the maximal cardinality of a set $S\ssq X$
which is {\em $\delta$-separated}\foot{Given any function
  $\rho:X\times X\to \R$, we call a set $S\ssq X$ $\delta$-separated
  with respect to $\rho$ if $\rho(x,y)\geq \delta$ for all $x\neq y\in
  S$. One should think of $\rho$ as a metric, but we will also use the
  same terminology in more general situations.} with respect to
$d^f_n$, then the {\em topological entropy} of $f$ can be defined by
\begin{equation}\label{e.entropy}
  \htop(f) \ = \ \lim_{\delta\to 0} \limsup_{n\to\infty} \log S_n(f,\delta)/n \ .
\end{equation}
This quantity measures the `chaoticity' of a dynamical system and is arguably
the most important topological invariant in ergodic theory.  If it is either
infinite or zero, however, then the complexity of a system has to be described
by different means. In the case of infinite entropy, mean dimension has been
established as a suitable substitute \cite{LindenstraussWeiss2000MeanDimension}.

If the entropy is zero, however, then the situation is less clear.  There exist
several alternative concepts to describe the complexity of a system in this
situation (see, for example, \cite{Misiurewicz1981,Smital1986,MisiurewiczSmital1988,
KolyadaSharkovsky1991,Ferenczi1997MeasureTheoreticComplexity,KatokThouvenot1997SlowEntropy,Ferenczi1999,
BlanchardHostMaas2000TopologicalComplexity,HasselblattKatok2002HandbookPrincipalStructures,
FerencziPark2007,HuangParkYe2007, ChengLi2010}), and different topological
invariants have been proposed for this purpose (\cite{Carvalho1997,
HasselblattKatok2002HandbookPrincipalStructures,HuangYe2009,DouHuangPark2011,
Marco2013,KongChen2014,FuhrmannGroegerJaeger2015AmorphicComplexity}).

Some of them have properties that may be considered as shortcomings, although
this partly depends on the viewpoint and the particular purpose one has in mind
(we briefly discuss this issue in Section~\ref{Discussion} below). In any case,
it is not always obvious which one should be considered best in a particular
situation, and in general there are still many gaps in the present state of
knowledge. At the same time, the issue has considerable relevance, since there
exist many system classes of both of theoretical and practical importance in
which the topological entropy is zero for structural reasons. Just to mention
some examples, these include regular Toeplitz flows
\cite{Downarowicz2005ToeplitzFlows}, circle homeomorphisms
\cite{katok/hasselblatt:1997}, interval exchange transformations \cite{Viana2006},
certain mathematical quasicrystals
\cite{Moody2000ModelSetsSurvey,BaakeLenzMoody2007Characterization},
quasiperiodically forced circle maps \cite{glendinning/jaeger/stark:2009} or
$\cC^{1+\alpha}$-surface diffeomorphisms with subexponential growth of periodic
orbits \cite{katok:1980}.

Maybe the most straightforward approach to the problem is to consider
subexponential, and in particular polynomial, growth rates instead of
exponential ones as in \eqref{e.entropy}. This leads to the notion of
{\em power entropy}\foot{defined in
  Section~\ref{Entropy_Definitions}} $\hpow$, which is also known under the
name of {\em polynomial word complexity} in the context of symbolic systems.
One aspect in which this quantity behaves quite differently from topological
entropy is the fact that it is very sensitive to transient behaviour.
For instance, the existence of a single wandering point\foot{We call
  $x\in X$ a {\em wandering point} of $f$ if there exists an open
  neighbourhood $U$ of $x$ such that $f^n(U)\cap U=\emptyset$ for all
  $n\in\N$.} of a homeomorphism $f$ implies $\hpow(f)\geq 1$ \cite{Labrousse2013}. In
particular, this means that the dynamically trivial Morse-Smale systems have
positive power entropy. A direct consequence is the non-existence of a variational
principle, which is another decisive difference to the standard notion
of topological entropy.

An alternative concept is {\em modified power entropy}
\cite{HasselblattKatok2002HandbookPrincipalStructures}.
In its definition, the Bowen-Dinaburg metrics are replaced by the corresponding
{\em Hamming metrics}. However, although this is less obvious to see, this
notion is equally sensitive to transient dynamics and therefore cannot satisfy a
variational principle either. We provide an example to demonstrate this
statement in Section~\ref{Counterexample}. Since this question has been left
open in the literature so far (see, for example, \cite[page 92]
{HasselblattKatok2002HandbookPrincipalStructures}), the communication of this
fact is one of the main motivations for this note.

The second issue we discuss here is the response of power entropy and
modified power entropy to the break of equicontinuity, which can be
observed during the transition from almost periodic
($=$equicontinuous) minimal systems to their almost 1-1 extensions. It
turns out that power entropy is suitable to detect this change in the
qualitative behaviour, whereas modified power entropy is not.  In this
context, we also introduce and discuss {\em amorphic complexity}. This
is a new topological invariant that equally measures the complexity
of zero entropy systems, but is based on an asymptotic rather than a
finite-time concept of separation \cite{FuhrmannGroegerJaeger2015AmorphicComplexity}.

Finally, we provide an example of an irregular Toeplitz flow with zero
modified power entropy in order to clarify some further aspects
of the preceding discussion.  \medskip

\noindent {\bf Acknowledgments.} Both authors have been supported by an
Emmy-Noether grant of the German Research Council (grant Oe 538/3-1).  T.J.~also
acknowledges support by a Heisenberg fellowship of the German Research Council
(grant OE 538/7-1). A substantial part of this work has been produced during a
stay of T.J. at the Max-Planck-Institute for Mathematics, Bonn, during the {\em
  `Dynamics and Numbers'}-activity in June 2014, and we would like to thank the
organisers for creating this opportunity.

\section{Some thoughts on slow entropies} \label{Discussion}

In the context of this discussion, we understand {\em `slow entropy'} in a broad
sense of a topological invariant that measures the complexity of dynamical
systems in the zero entropy regime. Thereby, our focus lies on $\Z$- and
$\N$-actions of low complexity. We note that similar concepts are also used for
the description of more general group actions, where the need for considering
alternative growth rates stems rather from the fast volume growth of the F\o
lner sequences than from the low complexity of the group action. However, we
will not go into any detail in this direction and refer to
\cite{KatokThouvenot1997SlowEntropy,Dong2014,KatokKatokHertz2014} for a
discussion and further references. In order to restrict the scope to some
degree, we concentrate on real-valued invariants and compact metric spaces. We
thus say a slow entropy is a function $h$ defined on the space of pairs $(X,f)$,
with $X$ a compact metric space and $f:X\to X$ continuous, which satisfies the
following requirements.

\begin{itemize}
\item $h$ is real-valued (including $\infty$) and non-negative;
\item $h$ takes the same value for topologically conjugate systems
  ({\em topological invariance});
\item If $g$ is a topological factor of $f$, then $h(g)\leq h(f)$
  ({\em monotonicity});
\item If $f$ has positive topological entropy, then $h(f)=\infty$
  ({\em zero entropy regime}).
\end{itemize}

Note that in fact topological invariance is a consequence of monotonicity.
Beyond these basic assumptions, however, it is not always clear what further
properties are desireable for a slow entropy, and which ones are rather not. The
reason behind is that this depends to a large extent on the purpose that such a
quantity should serve, and there are quite different and sometimes even
contradictory aims one could have in mind. We want to discuss this by means of
an example.

As it is well-known, one of the most important results about entropy is the
variational principle, which states that topological entropy equals the supremum
over its measure-theoretic counterparts with respect to all the invariant
probability measures of the system. It is one of the main tools in thermodynamic
formalism and explains the central role topological entropy plays in this
powerful machinery. As a consequence, topological entropy is also independent of
transient behaviour and determined by the dynamics on the set of recurrent
points only.  It is one possible aim for introducing a slow entropy to provide
similar tools for the study of zero-entropy systems. Most likely, however, this
will require at least some minimal amount of `chaoticity' in the system. In
contrast to this, an alternative task for a slow entropy would be to detect the
very onset of complicated dynamical behaviour. For example, one might want it to
detect the qualitative change in behaviour when going from equicontinuous
systems -- to which one would usually assign zero complexity -- to
non-equicontinuous systems, by taking positive values for the latter.  Now, this
would mean that, for instance, the slow entropy should give different values to
Sturmian subshifts and irrational rotations. However, a Sturmian subshift is
uniquely ergodic and measure-theoretically isomorphic to an irrational rotation,
so that this immediately contradicts a variational principle.

Hence, it seems obvious that there is not one single notion of slow entropy that
fulfills all the possible roles of a topological invariant in the zero entropy
regime at the same time. Certainly, this just reflects the great diversity of
zero entropy systems, which comprise many classes of quite different
complexity. The fact that not all of them can be adequately described with the
same concept is not too surprising. A more reasonable aim would be to identify a
whole array of useful invariants such that their union allow to cover the zero
entropy regime in a reasonable way and distinguish different degrees of
complexity. Yet, the present state of knowledge on the topic is still far from
this situation, and it will presumably need a lot of further fundamental
research in the area in order to get to that point.

The particular contribution of the present paper in this context is a
modest one. As mentioned, we concentrate mostly on modified power
entropy and clarify some of the mentioned aspects concerning this
particular notion. A short summary will be given in
Section~\ref{Conclusions}.

\section{Power entropy, modified power entropy and amorphic complexity}
\label{Entropy_Definitions}

As mentioned above, power entropy measures the polynomial growth rate
of orbits distinguishable by the Bowen-Dinaburg metrics $d^f_n$. In analogy to
\eqref{e.entropy}, it is defined as
\begin{equation*}
\hpow(f) \ = \ \lim_{\delta\to 0} \lim_{n\to\infty} \frac{\log S_n(f,\delta)}{\log n} \ ,
\end{equation*}
whenever the limits with respect to $n\to\infty$ exists. If this is not
the case, then one defines upper and lower power entropy $\hpowup$ and
$\hpowlow$ by taking the limit superior and limit inferior,
respectively. Note that due to the monotonicity in $\delta$, the
existence of the second limit is automatic.
We refer to \cite{Marco2013} for more information about this quantity. 

In the definition of {\em modified power entropy} (MPE), the
Bowen-Dinaburg metrics are replaced by the {\em Hamming metrics}
\begin{equation*}
\hat d_n^f(x,y) \ = \  \ntel\inergsum d(f^i(x),f^i(y)) \ .
\end{equation*}
If $\hat S_n(f,\delta)$ denotes the maximal cardinality of a set
$S\ssq X$ that is $\delta$-separated with respect to $\hat d^f_n$,
then the modified power entropy of $f$ is defined as 
\begin{equation*}
\hmod(f) \ = \ \lim_{\delta\to 0} \lim_{n\to\infty} \frac{\log \hat S_n(f,\delta)}{\log n} \ ,
\end{equation*}
provided the limit as $n\to\infty$ exists. If not, then one can again
define upper and lower versions $\hmodup$ and $\hmodlow$.  The fact
that $\hat d^f_n\leq d^f_n$ implies $\hat S_n(f,\delta)\leq
S_n(f,\delta)$ and hence $\hmod(f)\leq \hpow(f)$. We also note that
$\htop(f)>0$ implies $\hmod(f)=\infty$.\foot{This is well-known
  folklore, but for the convenience of the reader we provide a short
  direct proof in the next section.}

In both cases, the concept of separation that is used in the first step is one
in finite time: both metrics $d^f_n$ and $\hat d^f_n$ depend only on the first
$n$ iterates of the considered points. The limit for $n\to\infty$ is then taken
in a second step. However, since asymptotic notions like proximality, distality
or Li-Yorke pairs play a central role in topological dynamics, it seems natural
to also consider topological invariants that are directly based on an asymptotic
concept of separation. This is true for the following notion.

Given $x,y\in X$ and $\delta>0$, we let
\begin{equation*}
  M^{f}_{\delta,n}(x,y)\ = \ \#\left\{0\leq k<n\;|\;
    d(f^k(x),f^k(y))\geq\delta\right\} 
\end{equation*}
and 
\begin{equation*}
  \nu^f_\delta(x,y) \ = \ \varlimsup\limits_{n\to\infty} \frac{M^{f}_{\delta,n}(x,y)}{n} \ .
\end{equation*}
We say that $x$ and $y$ are \emph{$(f,\delta,\nu)$-separated} if
\[
   \nu^f_\delta(x,y) \ \geq \ \nu \ .
\]
Given $\nu>0$, we denote the maximal cardinality of a set $S\ssq X$
which is $\nu$-separated with respect to $\nu^f_\delta$ by
$S^*_\nu(f,\delta)$. Then, the \emph{amorphic complexity} of $f$ is defined as
\begin{equation*}
  \fsc(f) \ = \ \lim_{\delta\to 0} \lim_{\nu\to 0}\frac{\log S^*_\nu(f,\delta)}{-\log \nu} \ .
\end{equation*}
As before, this assumes that the limits with respect to $\nu$ exist.
Otherwise, it is again possible to define an upper and a lower
amorphic complexity.  Basic properties of this quantity, like
topological invariance, factor relations, power invariance and a product rule,
as well as the application to a number of example classes
are discussed in \cite{FuhrmannGroegerJaeger2015AmorphicComplexity}.
Somewhat surprisingly, amorphic complexity turns out be very well
applicable and accessible to explicit computations in various system
classes like regular Toeplitz flows, Sturmian shifts and Denjoy type
circle homeomorphisms or cut and project quasicrystals. The reason
behind is the fact that the asymptotic nature of the employed
separation concept allows to obtain bounds on the separation numbers
$S^*_\nu(f,\delta)$ by applying suitable ergodic theorems. We refer to
\cite{FuhrmannGroegerJaeger2015AmorphicComplexity} again for details.

The main reason for treating amorphic complexity here is to complete the
discussion in \cite[Section 3.7]{FuhrmannGroegerJaeger2015AmorphicComplexity} by
showing that there are no direct relations, in terms of an inequality, between
amorphic complexity and the other two notions.  Thereby, for one of the
directions, we will have to rely on the same example as for the non-existence of
a variational principle for modified power entropy. Hence, we come back to this
issue at the end of the next section.

In all of the above, we have considered polynomial growth rates, which
turn out to be the appropriate scale for many important example
classes. In general, however, it is certainly possible to take into
account more or less arbitrary growth rates. We say $a:\R_+\times \N
\to\R_+$ is a {\em scale function} if $a$ is strictly increasing in
both arguments.  Then, in analogy to the power entropy, the {\em upper
  $a$-entropy} of $f$ is defined as
\[
\overline h_a(f) \ = \ \sup_{\delta>0} \sup\left\{s>0 \ \left| \ \varlimsup_{n\to \infty}
    \frac{S_n(f,\delta)}{a(s,n)} > 0 \right.\right\} \ ,
\]
and the lower one accordingly.
In order to obtain good properties, one usually assumes that the scale
functions are $O$-regularly varying, that is, $\varlimsup_{n\to\infty}
\frac{a(s,mn)}{a(s,n)} <\infty$ for all $m\in\N$. Since this
definition allows to capture any rates of asymptotic growth in the
subexponential regime, one of the most important distinctions on the
qualitative level is whether $\sup_{n\in\N} S_n(f,\delta)$ is bounded
for all $\delta>0$, or infinite for all sufficiently small $\delta$.
We will mainly focus on this aspect in our discussion of almost 1-1
extensions in Section~\ref{Extensions}. Of course, all these comments
on the use of different scale functions equally apply to modified
power entropy and amorphic complexity. For the latter, scale functions
need to have the separation frequency as the second argument, so in
this case $a$ is a positive real-valued function on $\R^+\times (0,1]$ and
$O$-regularly varying means that $\varlimsup_{\nu\to 0}
\frac{a(s,c\nu)}{a(s,\nu)}$ is finite for all $c>0$.

\section{Modified power entropy and topological entropy}

In many cases, modified power entropy is strictly smaller than power
entropy (see e.g. \cite{HasselblattKatok2002HandbookPrincipalStructures}).
However, on an exponential scale (that is, using the scale function
$a(s,n)=\exp(sn)$), this difference disappears. Since we do not know
an appropriate reference, we include a proof of this well-known
result.

\begin{lemma} \label{l.entropy}
  Suppose $X$ is a compact metric space and $f:X\to X$ is continuous. Then
\[
\lim_{\delta\to 0}\varliminf_{n\to\infty} \frac{\log
  \hat S_n(f,\delta)}{n} \ = \ \lim_{\delta\to 0}\varlimsup_{n\to\infty}
\frac{\log \hat S_n(f,\delta)}{n} \ = \ \htop(f) \ .
\]
In particular, $\htop(f)>0$ implies $\hmod(f)=\infty$.
\end{lemma}
\proof The $\leq$-inequalities are obvious. Further, it is well-known
that $$\htop(f) = \lim_{\eps\to 0}\varliminf_{n\to\infty} \frac{\log
  S_n(f,\eps)}{n}$$ (e.g.~\cite{katok/hasselblatt:1997}). Therefore,
it suffices to show that
\begin{equation}
  \label{eq:1}
  \lim_{\delta\to 0}\varliminf_{n\to\infty} \frac{\log
  \hat S_n(f,\delta)}{n} \ \geq \ \lim_{\eps\to 0}\varliminf_{n\to\infty} \frac{\log
  S_n(f,\eps)}{n} \ . 
\end{equation}
To that end, fix $\eps>0$ and $\alpha>0$ and choose
$\delta\in(0,\eps\alpha/2)$. Further, let $U_1\ld U_K$ be a finite
partition of $X$ into sets of diameter $< \eps$. 

Given $n\in\N$, let $N=\hat S_n(f,\delta)$ and choose a partition of $X$
into sets $P_1\ld P_N$ with the property that $\hat d^f_n(x,y)\leq \delta$
for all $x,y\in P_j,\ j=1\ld N$. From each of the $P_j$, we select one
point $x_j\in P_j$. (Note that all of the $P_j$ are non-empty due to
the definition of $N$.) Then, given $x\in P_j$, we define $\omega(x)\in
\{0\ld K\}^n$ by
\[
\omega_i(x) \ = \left\{
  \begin{array}{cl}
    0 & \textrm{if } d(f^i(x),f^i(x_j)) \  < \ \eps/2 \ , \\ \ \\
    k & \textrm{if } d(f^i(x),f^i(x_j)) \ \geq \ \eps/2 \textrm{ and } x\in U_k  \ .
  \end{array} \right. 
\]
Note that if $x,y\in P_j$ and $\omega(x)=\omega(y)$, then
$d^f_n(x,y)<\eps$. Hence, the maximal cardinality of a subset of $P_j$
which is $\eps$-separated with respect to $d^f_n$ is at most
$\#\{\omega(x)\mid x\in P_j\}$. Moreover, we have that for each $x\in
P_j$
\[
\#\{0\leq i\leq n-1 \mid \omega_i(x)\neq 0\} \ \leq \ \alpha n \ ,
\]
since otherwise $\hat d^f_n(x,x_j) \geq \alpha n \eps/2 \ > \ \delta$,
contradicting the choice of the $P_j$. Hence, using that
$\twovector{n}{\lfloor n\alpha\rfloor} \leq \exp(-\alpha\log(\alpha)
n)$ for sufficiently small $\alpha$, we obtain that 
\begin{eqnarray*}
  \#\{\omega(x)\mid x\in P_j\} \ \leq \ \twovector{n}{\lfloor n\alpha\rfloor} 
  \cdot K^{\lfloor n\alpha\rfloor} \ \leq \ \exp\left( \alpha(\log(K)-\log(\alpha))n\right) \ . 
\end{eqnarray*}
Altogether, this yields that 
\[
S_n(f,\eps) \ \leq \ \hat S_n(f,\delta) \cdot \exp\left( \alpha(\log(K)-\log(\alpha))n\right) \ .
\]
Since $\lim_{\alpha\to 0} \alpha(\log(K)-\log(\alpha)) =0$ and
$\alpha>0$ was arbitrary, this proves (\ref{eq:1}).  \qed\medskip

\section{A counterexample to the existence of a variational principle
  for MPE} \label{Counterexample}

Let $I=[0,1]$ and $\kreis=\R/\Z$. The main aim in this section is the
construction of an example of the following type.
\begin{theorem} \label{t.counterexample}
  There exists a skew product map of the form
\begin{equation*}
  f:I\times\kreis\to
  I\times\kreis\quad ,\quad  f(x,y)=(\tau(x),y+\beta(x)+\rho) \ ,
\end{equation*}
where 
\begin{itemize}
\item $\rho\in\R\smin \Q$,
\item $\tau$ is a diffeomorphism of $I$ with exactly two
fixed points at $0$ and $1$,
\item $\beta:I\to \kreis$ is a differentiable function
with $\beta_{|\{0\}\cup B_\eps(1)}=0$ for some $\eps>0$,
\end{itemize}
such that $f$ satisfies
\[
\hmodlow(f)\geq 1/2 \ .
\]
\end{theorem}
Before we turn to the proof, we first draw the following conclusion.
\begin{corollary} \label{c.no_VP} There is no real-valued isomorphism
  invariant of measure-preserving dynamical systems that satisfies a
  variational principle with modified power entropy.
\end{corollary}
There is, of course, a standard measure-theoretic analogue of modified
power entropy, introduced in \cite{Ferenczi1997MeasureTheoreticComplexity,KatokThouvenot1997SlowEntropy}
(see also \cite{HasselblattKatok2002HandbookPrincipalStructures}), which is
bounded above by topological modified power entropy. However, since the only
structural property that is needed is the invariance under isomorphisms, we do
not need to state any detail here and omit these for the sake of brevity.

We also note that we understand `variational principle' in the sense that the
topological quantity equals the supremum over all measure-theoretic ones,
where the supremum is taken over all {\em invariant} measures. For the standard
notion of entropy it suffices to consider only ergodic measures due to the
linearity of measure-theoretic entropy, but in general this can make a big
difference (see also \cite[Section 4.4b and page 81]
{HasselblattKatok2002HandbookPrincipalStructures}).

\proof[Proof of Corollary \ref{c.no_VP}.]  Suppose that $h^*$ is a real-valued
function of pairs $(f,\mu)$, where $f$ is a continuous map on some (compact)
metric space.  Suppose for a contradiction that $\hmod(f)=\sup_{\mu}
h^*(f,\mu)$, where the supremum is taken over all invariant measures $\mu$ of
$f$.

In the example in Theorem~\ref{t.counterexample}, there exist exactly two
ergodic measures $\mu_0$ and $\mu_1$, which are the one-dimensional Lebesgue
measures on the two circles $\cT_0=\{0\}\times \kreis$ and
$\cT_1=\{1\}\times\kreis$, and the restriction of $f$ to these circles is just
the rotation by $\rho$.  Any invariant measure $\mu$ is a convex combination of
$\mu_0$ and $\mu_1$ and obviously isomorphic to $(f_{|\cT_0\cup
  \cT_1},\mu_{|\cT_0\cup\cT_1})$. However, as $f$ restricted to $\cT_0\cup\cT_1$
is an isometry we have $\hmod(f_{|\cT_0\cup\cT_1})=0$ and hence $h^*(f,\mu)=0$
by the assumed variational principle. This means $\sup_{\mu}h^*(f,\mu)=0$,
whereas $\hmodlow(f)\geq 1/2$, which yields the required contradiction.
\qed\medskip

We also note that in the situation of Theorem~\ref{t.counterexample} we have
$\hmodlow(f) > \hmodlow(f_{|\Omega(f)})$, where $\Omega(f)$ denotes the set of
non-wandering points of $f$. This shows that modified power entropy is sensitive
to transient dynamics. Since all invariant measures are supported on the
non-wandering set, this is equally not compatible with a variational
principle. We turn to the construction of the example.

\proof[Proof of Theorem~\ref{t.counterexample}.]

We first construct the diffeomorphism $\tau :I\to I$. To that end,
let $I_n=[2^{-n},3\cdot 2^{-(n+1)}]$ and $I_n'=[2^{-n},5\cdot
2^{-(n+2)}]$ where $n\geq 3$. Then, we choose a $\cC^1$-function
$\alpha:I\to I$ with the following properties:
\begin{itemize}
\item[(i)] $|\alpha'(x)| < 1$ for all $x\in I$;
\item[(ii)] $\alpha(0)=\alpha(1)=0$;
\item[(iii)] $\alpha(x)>0$ for all
$x\in(0,1)$;
\item[(iv)] $\alpha_{|I_n} = 2^{-(3n+4)}$;
\end{itemize}
Further, we let $\beta(x)=x$ if $x\in[0,7/8]$ and extend this differentiable to
all of $I$ in such a way that $\beta_{|B_\eps(1)}=0$ for some $\eps>0$. Note that
since $d(I_n,I_{n+1})= 2^{-(n+2)}$, condition (iv) does not contradict the
differentiability of $\alpha$. Due to (i) and (ii) the map $\tau:I\to I,\
x\mapsto x+\alpha(x)$ is a $\cC^1$-diffeomorphism of $I$ with unique fixed
points $0$ and $1$. Moreover, due to (iii) we have $\nLim\tau^n(x)=1$ for all
$x\in (0,1]$.

In order to prove $\hmodlow(f)>0$, fix $n\geq 3$ and choose $x_1^n <
x_2^n < \ldots < x_{2^n}^n \in I_n'$ with
$x_{i+1}^n-x_i^n=2^{-(2n+2)}$. By (iv) and the choice of the
intervals $I_n'$ and $I_n$, we have $\tau^k(x^n_j)\in I_n$ for all
$j=1\ld 2^n$ and $k=0\ld 2^{2n+2}$. Since $\alpha$ is constant on
$I_n$, this means that the points $x^n_j$ remain at equal distance for
the first $2^{2n+2}$ iterations. If we consider the $2^n$ points
$(x^n_j,0)\in I\times\kreis$, then for $l,m=1\ld 2^n$ the vertical distance
after $n$ steps is
\[
d\left(\pi_2\circ f^k(x^n_l,0),\pi_2\circ f^k(x^n_m,0)\right) \ = \ d'\left(k\cdot
(l-m)\cdot 2^{-(2n+2)},0\right) \  .
\]
Here $d'$ denotes the canonical distance on $\kreis$. An easy
computation yields for $l\neq m$
\[
\hat d^f_{2^{2n+2}}((x^n_l,0),(x^n_m,0)) \ \geq \ \frac{1}{2^{2n+2}}
\sum_{k=0}^{2^{2n+2}-1}d'\left(k\cdot (l-m)\cdot 2^{-(2n+2)},0\right) \
= \ \viertel \ .
\]
This means that the set $\{x^n_1\ld x^n_{2^n}\}$ is
\viertel-separated with respect to $\hat d^f_{2^{2n+2}}$. We thus
obtain $\hat S_n(f,1/4)\geq 2^n$ and hence
\[
\hmodlow(f) \ \geq \ \varliminf_{n\to\infty} \frac{\log
  \hat S_n(f,\delta)}{\log n} \ \geq \ \kLim
\frac{\log(2^k)}{\log(2^{2k+2})} \ = \halb \ .\qedhere
\]
\medskip

In order to conclude this section, we want to discuss why the above
example also shows that there is no direct relation, in terms of an
inequality, between modified power entropy and amorphic complexity.
To that end, let us first look at some trivial examples. 

Since Morse-Smale systems have a finite set of fixed or periodic points and these
attract all other orbits, it is easy to see that they have zero amorphic
complexity. This shows that one may have $\fsc(f)<\hpow(f)$. On the other hand,
consider $f:\torus\to\torus,\ (x,y)\mapsto (x,x+y)$. Then any two points with
different $x$-coordinate rotate with different speed in the vertical direction,
and it is therefore easy to see that they are $\nu$-separated with respect to
$\nu^f_{\delta}$ if $\nu,\delta>0$ are chosen sufficiently small.  Therefore,
$S^*_\nu(f,\delta)=\infty$ and thus $\fsc(f)=\infty$. At the same time, it is
easy to check that $\hpow(f)\leq 1$ (see \cite[Section
3.7]{FuhrmannGroegerJaeger2015AmorphicComplexity}). Hence, we may have
$\fsc(f)>\hpow(f)$ (and thus also $\fsc(f)>\hmod(f)$). The only remaining
direction is therefore to show that $\hmod(f)>\fsc(f)$ is possible as
well. However, we claim that this is the case in the example constructed
above. In order to see this, the following basic observation is helpful.
\begin{lemma}[{\cite[Lemma 3.11]{FuhrmannGroegerJaeger2015AmorphicComplexity}}]
  Suppose that $X$ is a compact metric space, $f$ is a continuous map and $A\ssq
  X$ is a forward invariant subset such that for all $x\in X\smin A$ there exists
  $y_x\in A$ such that $\nLim d(f^n(x),f^n(y_x))=0$. Then
  $\fsc(f)=\fsc(f_{|A})$.\label{l.unique_target}
\end{lemma}
In the above situation, we say that $f$ has the {\em unique target property}
with respect to $A$. If $f$ is the example constructed in the proof of
Theorem~\ref{t.counterexample}, then this assumption is satisfied for
$A=\Omega(f)$. This is a direct consequence of the fact that
$\beta_{|B_\eps(1)}=0$ (and the reason for including this condition, which has
not been used otherwise). Note here that all orbits outside of $\cT_0$ converge
to the circle $\cT_1$ upon forward iterations, and once they enter
$B_\eps(1)\times\kreis$ the rotation in the second coordinate is always equal to
$\rho$. For this reason, all these orbits have a unique `target orbit' in
$\cT_1$. Since $f_{|\cT_0\cup\cT_1}$ is an isometry and therefore has amorphic
complexity zero, the above Lemma~\ref{l.unique_target} yields
$\fsc(f)=0<\hmodlow(f)$.

It remains to point out that since the example constructed above has the unique
target property with respect to the non-wandering set, the transient dynamics
causing the positive modified power entropy should still be considered as rather
`tame'. Amorphic complexity is equally sensitive to transient dynamics, but
these have to `mix up' orbits arbitrarily close to the non-wandering
set. An example similar to the one above is given in \cite[Section
3.5]{FuhrmannGroegerJaeger2015AmorphicComplexity}.

\section{Modified power entropy of regular almost 1-1 extensions}
\label{Extensions}

Given two compact metric spaces $(X,d), (\Xi,\rho)$ and two continuous maps
$f:X\to X$, $\tau:\Xi\to\Xi$, we say $(X,f)$ is a {\em (topological) extension} of
$(\Xi,\tau)$ if there exists a continuous onto map $h:X\to\Xi$ such that $h\circ
f=\tau\circ h$. In this situation, $h$ is called a {\em factor map} or {\em
  semi-conjugacy} from $f$ to $\tau$ and $(\Xi,\tau)$ is called a {\em
  (topological)} factor of $(X,f)$.  For the sake of brevity, we will sometimes
omit the spaces and say $\tau$ is a factor of $f$. An extension is called {\em
  almost 1-1} if the set $\Omega=\{\xi\in\Xi\mid \#h^{-1}(\xi)=1\}$ is generic
in the sense of Baire (that is, a residual set). Note that if $f$ and $\tau$ are
invertible, then the set $\Omega$ is $\tau$-invariant. Moreover, if in addition
$\tau$ is minimal, then it suffices to require that there exist a single $\xi$
with $\# h^{-1}(\xi)=1$.

From now on, we assume for the remainder of this section that $f$ and $\tau$ are
invertible and $\tau$ is minimal. Further, we suppose that $\tau$ is {\em almost
  periodic} (that is, equicontinuous). In this case, there exists a unique
$\tau$-invariant probability measure $\mu$ on $\Xi$, which is necessarily
ergodic (unique ergodicity). We say that the extension $(X,f)$ is {\em regular}
if $\mu(\Omega)=1$ and {\em irregular} if $\mu(\Omega)=0$. Note that by
invariance of $\Omega$ and ergodicity of $\mu$, one of the two always holds.  We
refer to \cite{auslander1988minimal} for a comprehensive exposition.

A regular almost 1-1 extension of an equicontinuous minimal system is always
uniquely ergodic and isomorphic to its factor. For this reason, the topological
entropy is zero in this case.  However, if the extension is not everywhere 1-1
(that is, there exists $\xi\in\Xi$ such that $\# h^{-1}(\xi)>1$), then $f$
cannot be equicontinuous.\foot{From now on, whenever speaking of extensions we
  will assume implicitly that the factor map is not injective.} Thus, there is a
break of equicontinuity when going from equicontinuous minimal systems to their
almost 1-1 extensions, but at the same time this does not lead beyond the regime
of zero entropy. It is therefore a natural question to ask how a topological
invariant for low-complexity systems behaves during this bifurcation. In
particular, it is one possible task for such a slow entropy to detect this
change in the qualitative behaviour. We will discuss a positive result in this
direction for amorphic complexity and power entropy further below. Modified
power entropy, however, does not respond to this transition.

\begin{theorem}\label{t.almost_sure_extensions}
  Suppose $f:X\to X$ is a regular almost 1-1 extension of a minimal
  equicontinuous homeomorphism $\tau:\Xi\to\Xi$. Then $\sup_{n\in\N} \hat
  S_n(f,\delta)<\infty$ for all $\delta>0$ and in particular $\hmod(f)=0$.
\end{theorem}
For the proof, the following statement will be useful. 
\begin{lemma}[{\cite[Lemma 2.6]{FuhrmannGroegerJaeger2015AmorphicComplexity}}]\label{l.eta}
  Let $h:X\to\Xi$ be the factor map of an almost 1-1 extension and define
  \[ 
   E_\delta \ = \ \{\xi\in\Xi\mid
  \diam(h^{-1}(\xi))\geq \delta\} \ .
  \]
 Then for all $\delta>0$ and $\eps>0$ there
  exists $\eta_\delta(\eps)>0$ such that if $x,y\in X$ satisfy $d(x,y)\geq\delta$
  and $\rho(h(x),h(y))<\eta_\delta(\eps)$, then $h(x)$ and $h(y)$ are both
  contained in $B_\eps(E_\delta)$.
\end{lemma}
\proof[Proof of Theorem~\ref{t.almost_sure_extensions}.] By going over to an
equivalent metric, we may assume without loss of generality that $\tau$ is an
isometry. Fix $\delta>0$ and choose
\begin{equation}\label{e.nu-choice}
\nu \ < \ \frac{\delta}{2\diam(X)-\delta} \ .
\end{equation}
Then, choose $\eps>0$ such that $\mu\left(A\right)<\nu$ where
$A=\overline{B_{\eps}(E_{\delta/2})}$. Let $\eta=\eta_{\delta/2}(\eps)$
  be as in Lemma~\ref{l.eta}. Due to the Uniform Ergodic Theorem we
  can find $M\in\N$ such that for all $n\geq M$ and $\xi\in\Xi$ we
  have
\begin{equation*}
  \ntel \inergsum \ind_A\circ 
  \tau^i(\xi) \ < \ \nu \ .
\end{equation*}
Therefore, given two points $x,y\in X$ with $\rho(h(x),h(y))<\eta$ and
$n\geq M$, the fact that $\tau$ is an isometry together with
Lemma~\ref{l.eta} implies 
\begin{eqnarray*}
  \hat d^f_n(x,y) & \leq & \frac{\diam(X)}{n} \cdot \left(\inergsum
    \ind_A\circ \tau^i(h(x))\right) +  \frac{\delta}{2n}\cdot
  \left(n- \inergsum \ind_A\circ \tau^i(h(x))\right) \\
  & \leq &\ \diam(X)\cdot \nu + \frac{\delta}{2}\cdot (1-\nu) \ \stackrel{(\ref{e.nu-choice})}{<}
     \ \delta \ .
\end{eqnarray*}
Thus, independent of $n\geq M$, two points $x,y$ can only be $\delta$-separated
with respect to $\hat d^f_n$ if $h(x)$ and $h(y)$ have distance greater than
$\eta$ in $\Xi$. Hence, any set $S\ssq X$ which is $\delta$-separated with
respect to $\hat d^f_n$ projects to a set which is $\eta$-separated with respect
to the metric in $\Xi$. Since $\Xi$ is compact, there exists an upper bound
$K(\eta)$ on the maximal cardinality of an $\eta$-separated set. We obtain $\hat
S_n(f,\delta)<K(\eta)$ for all $n\geq M$ and consequently $\hmod(f)=0$ as
claimed.\qed\medskip

Theorem~\ref{t.almost_sure_extensions} is in contrast to the following result in
\cite{FuhrmannGroegerJaeger2015AmorphicComplexity}, which shows that the
asymptotic separation numbers involved in the definition of amorphic complexity
are sensitive to the break of equicontinuity in the above situation.

\begin{theorem}[{\cite[Theorem 2.10]{FuhrmannGroegerJaeger2015AmorphicComplexity}}]
  Let $f:X\to X$ be a minimal almost 1-1 extension of an equicontinuous homeomorphism
  $\tau:\Xi\to\Xi$. Then there is $\delta>0$ such that
  $\sup_{\nu>0}  S_\nu^*(f,\delta)=\infty$.
\end{theorem}
We note that this result does not guarantee polynomial growth rates, but as
discussed in Section~\ref{Entropy_Definitions} one can obtain positive amorphic
complexity with respect to a suitably chosen scale function. In order to obtain
a similar result for power entropy, it suffices to use the following elementary
observation.
\begin{lemma}
  If $\sup_{\nu>0} S_\nu^*(f,\delta)=\infty$, then $\sup_{n\in\N}
  S_n(f,\delta)=\infty$.
\end{lemma}
\proof Suppose that $S$ is a set which is $\nu$-separated with respect to
$\nu^f_\delta$. Then there exists $n>0$ such that for all $x\neq y\in S$ we have
$\# \{i=0\ld n-1\mid d(f^i(x),f^i(y))\geq \delta\}/n \geq \nu/2 >0$. This
immediately implies $d^f_n(x,y)\geq \delta$ for all $x\neq y\in S$. Hence, $S$
is a $\delta$-separated set with respect to $d^f_n$ and therefore
$S_n(f,\delta)\geq\#S$. The statement now follows easily. \qed\medskip

We note that the following direct consequence is also contained in a more
general result by Blanchard, Host and Maass \cite[Proposition
2.2]{BlanchardHostMaas2000TopologicalComplexity}.

\begin{corollary}
  Suppose $f:X\to X$ is a minimal almost 1-1 extension of an equicontinuous
  homeomorphism $\tau:\Xi\to\Xi$. Then there exists $\delta>0$ such that
  $\sup_{n\in\N}  S_n(f,\delta)=\infty$.
\end{corollary}

\section{An example of an irregular Toeplitz flow of low complexity}
\label{Toeplitz}

   Since modified power entropy does not detect the difference between
   equicontinuous minimal systems and their regular almost 1-1 extensions, one
   could hope that instead it responds to the transition from regular to
   irregular almost 1-1 extensions. The aim of this section is to demonstrate
   that this is not the case either.  To that end, we construct an example of an
   irregular Toeplitz sequence, leading to an irregular almost 1-1 extension $f$
   of the corresponding odometer, which has modified power entropy zero and even
   bounded separation numbers $\hat S_n(f,\delta)$ for all $\delta>0$.  We
   assume some acquaintance with the theory of odometers and Toeplitz flows and
   refer to the excellent survey \cite{Downarowicz2005ToeplitzFlows} or classical papers by
   Jacobs and Keane \cite{JacobsKeane1969}, Eberlein \cite{Eberlein1971} and
   Williams \cite{Williams1984} for the relevant details.

We let $\Sigma=\{0,1\}^\Z$ and equip it with the metric
\begin{equation*}
  d(\omega,\tilde\omega) \ = \ \sum_{\omega_k\neq\tilde\omega_k }
  2^{-|k|} \ ,
\end{equation*}
where $\omega,\tilde\omega\in\Sigma$ and the index $k$ in the sum runs over all
of $\Z$, to make it a compact metric space. By $\sigma$ we denote the left shift
on $\Sigma$. Given $\omega\in\Sigma$, we let $\Sigma_\omega$ be the shift orbit
closure of $\omega$.

\begin{theorem} \label{t.ToeplitzExample} There exists an irregular
  Toeplitz sequence $\omega$ such that the corresponding Toeplitz flow
  $(\Sigma_\omega,\sigma)$ satisfies $\sup_{n\in\N}\hat
  S_n(\sigma_{|\Sigma_\omega},\delta) < \infty$ for all $\delta>0$,
  and in particular $\hmod(\sigma_{|\Sigma_\omega})=0$.
\end{theorem}
Note that conversely irregular Toeplitz flows may have positive
entropy \cite{BulatekKwiatkowski1992}, in which case the modified
power entropy is infinite by Lemma~\ref{l.entropy}.

Before we turn to the proof, we first need to address some technical
issues.  Given $\omega,\tilde\omega\in\Sigma$ and $n\in\N$, we let
\[
D_n(\omega,\tilde\omega) \ = \ \frac{1}{2n+1} \sum_{i=-n}^n |\omega_i-\tilde\omega_i| \ .
\]
Further, we denote by $R_n(\Sigma_\omega,\delta)$ the largest cardinality of a
set $R\ssq\Sigma_\omega$ which is $\delta$-separated with respect to $D_n$. The
following statements allow to relate $R_n(\Sigma_\omega,\delta)$ to $\hat
S_n(\sigma_{|\Sigma_\omega},\delta)$.
\begin{lemma}
  Let $n\in\N$.
  We have $\hat d^\sigma_n(\omega,\tilde\omega) \leq 9D_{2n}(\omega,\tilde\omega) + 2^{-(n-1)}$
  for all $\omega,\tilde\omega\in \Sigma$. In particular, if
  $\hat d^\sigma_n(\omega,\tilde\omega)\geq \delta$ and $2^{-(n-1)}\leq \delta/2$, then
  $D_{2n}(\omega,\tilde\omega)\geq \delta/18$ and hence
  $R_{2n}(\Sigma_\omega,\delta/18) \geq \hat S(\sigma_{|\Sigma_\omega},\delta)$.
\end{lemma}
\proof We have
\begin{eqnarray*}
  \hat d^\sigma_n(\omega,\tilde\omega) & = & \ntel\inergsum d(\sigma^i(\omega),\sigma^i(\tilde\omega)) \
  = \    \ntel \inergsum \sum_{\omega_{k+i}\neq\tilde\omega_{k+i}} 2^{-|k|} \\ 
  & = & \ntel  \sum_{\omega_k\neq\tilde\omega_k} \inergsum 2^{-|k-i|}  \ \leq \
  \ntel\left( \sum_{|k|\leq 2n: \omega_k\neq\tilde\omega_k} 3 \ + 
\sum_{|k|>2n:\omega_k\neq\tilde\omega_k} 2^{-(|k|-n)} \right) \\
  & \leq & \frac{3(2n+1)}{n}D_{2n}(\omega,\tilde\omega) + 2^{-(n-1)} 
\ \leq \ 9D_{2n}(\omega,\tilde\omega) + 2^{-(n-1)} \ .
\end{eqnarray*}
\qed\medskip

\begin{corollary} \label{c.symbolic_separation_oount} 
  If $\sup_{n\in\N} R_n(\Sigma_\omega,\delta)<\infty$ for all $\delta>0$,
  then $\sup_{n\in\N} \hat
  S_n(\sigma_{|\Sigma_\omega},\delta)<\infty$ for all $\delta>0$.
\end{corollary}

Thus, it suffices to consider the pseudometrics $D_n$ in the proof of
Theorem~\ref{t.ToeplitzExample}, which are easier to handle than the
Hamming metrics in this context. 

For the particular case of Toeplitz flows, there is a further simplification,
which is due to the fact that for a Toeplitz sequence $\omega$ the space
$\Sigma_\omega$ consists precisely of those sequences which have exactly the
same subwords as $\omega$. This leads to the following elementary
observations. In order to specify finite subwords of $\omega\in\Sigma$, we let
$\omega^{m,n}_j=\omega_{m+j}$ for $j=-n\ld n$, so that $\omega^{m,n}$ is the
subword of $\omega$ with length $2n+1$ and center position $m$.  In order to
count the number of mismatches between two subwords of the same length, we let
\[
D_n(\omega^{m,n},\omega^{m',n})\ =\ D_n(\sigma^m(\omega),\sigma^{m'}(\omega))
\ = \ \frac{1}{2n+1} \sum_{j=-n}^n |\omega^{m,n}_j-\omega^{m',n}_j| \ .
\]
Further, we denote by $\tilde R_n(\omega,\delta)$ the largest
cardinality of a family of subwords of $\omega$ of length $2n+1$ which
are $\delta$-separated with respect to $D_n$. 
\begin{corollary}
  \label{l.finite_words} If $\omega\in\Sigma$ is a Toeplitz sequence, then we
  have $\tilde R_n(\omega,\delta) = R_n(\Sigma_\omega,\delta)$ for all $n\in\N$,
  $\delta>0$. In particular, if for each $\delta>0$ we have 
  $\sup_{n\in\N} \tilde R_n(\omega,\delta)<\infty$, then
  $\sup_{n\in\N}\hat S_n(\sigma_{|\Sigma_\omega},\delta)<\infty$ for all $\delta>0$.
\end{corollary}
We can now turn to the

\proof[Proof of Theorem~\ref{t.ToeplitzExample}.] Our construction is a classical
one which has been used in similar form by many authors \cite{Oxtoby1952,
Williams1984,BulatekKwiatkowski1990,BulatekKwiatkowski1992,Downarowicz2005ToeplitzFlows}.
The difficulty lies in controlling the separation numbers with respect to the
Hamming metrics.

We first fix $a_1\in\N$ and a sequence $(b_n)_{n\in\N}$ of integers $\geq 2$,
specified further below, and let $a_{n+1}=2b_na_n$ for all $n\geq 1$. Further,
we let
\[
A_n \ = \ \{-a_n\ld a_n\} + a_{n+1}\Z \quad , \quad B_n \ = \ \bigcup_{i=1}^n
A_n \quad \textrm{and} \quad C_n\ = \ B_n\smin B_{n-1} \ .
\]
Intervals of the form $\{-a_n\ld a_n\}+\ell a_{n+1}$ with $\ell\in\Z$
will be called {\em $n$-blocks}. Note that since the $a_n$ converge to
$\infty$, we have $\biguplus_{n\in\N} C_n=\Z$. Now, we can define a
Toeplitz sequence $\omega\in\Sigma$ by
\[
\omega_k \ = \ \left\{
  \begin{array}{ccl}
    0 & \textrm{if} & k\in C_n \textrm{ with } n \textrm{ odd};\\ \ \\
    1 & \textrm{if} & k\in C_n \textrm{ with } n \textrm{ even}.\\
  \end{array} \right. \ 
\]
If $k\in C_n$, then by definition $k$ is an $a_{n+1}$-periodic
position, that is, $\omega_{k+\ell a_{n+1}}=\omega_k$ for all
$\ell\in\Z$. By construction, all positions $k\in\Z$ are periodic for
some $a_n$ in this sense, so that by definition $\omega$ is a Toeplitz
sequence with periodic structure $(a_n)_{n\in\N}$.  Further, if
$k\in C_{n+1}$, then there exists $\ell\in\Z$ with $k+\ell a_{n+1}\in
C_{n+2}\smin C_{n+1}$. Hence, we have $\omega_{k+\ell a_{n+1}}\neq
\omega_k$, so that $k$ is not an $a_{n+1}$-periodic position.
Therefore, we obtain that the set $\Per(\omega,a_{n+1})$ of
$a_{n+1}$-periodic positions equals $C_n$.  Since $C_n=\bigcup_{i=1}^n
A_i$, we obtain that this set has density
\[
\cD(a_{n+1}) \ \leq \ \sum_{i=1}^{n} \frac{2a_i}{a_{i+1}} \ = \
\sum_{i=1}^n \frac{1}{b_i} \ .
\]
If we choose the $b_i$'s such that $\sum_{i=1}^\infty b_i^{-1}<1$, then $\nLim
\cD(a_n) < 1$. This means, by definition, that the Toeplitz sequence $\omega$
with periodic structure \nfolge{a_n} is {\em irregular} and thus
$(\Sigma_\omega,\sigma)$ is an irregular almost 1-1 extension of a corresponding
odometer (see \cite{Downarowicz2005ToeplitzFlows}).

In order to show that $\sup_{n\in\N} \tilde R_n(\omega,\delta)
<\infty$ for all $\delta$, we fix $\delta>0$ and $j_0\in\N$ with
$2^{-j_0}< \delta/4$. Since $\tilde R_n(\omega,\delta)\leq \tilde
R_{n'}(\omega,\delta/2)$ if $n\leq n'\leq 2n$, it suffices to show
that
\begin{equation}\label{e.toeplitz_proof}
\sup_{j\in\N} \tilde R_{2^j}(\omega,\delta) \ < \ \infty
\end{equation}
for all $\delta>0$. To that end, choose $s\in\N$ with $\sum_{i=s+1}^\infty
b_i^{-1}<\delta/8$, so that the asymptotic density of $B_n\smin B_s$ is smaller
than $\delta/4$ for all $n>s$.  Given $N=2^j$ with $j\geq j_0$, we will define a
partition $\{A_{p,q}^\iota\}_{p,q}^\iota$ of $\Z$ with the property that if
$m,m'$ belong to the same element of the partition, then $\omega^{m,N}$ and
$\omega^{m',N}$ cannot be $\delta$-separated with respect to $D_n$. This implies
immediately that $\tilde R(\omega,\delta)$ does not exceed the number of
partition elements.  Since the latter will be independent of $N$, this will
prove (\ref{e.toeplitz_proof}).

The partition elements $\cA_{p,q}^\iota$ depend on three parameters. The parameter
$p$ describes the position of the subwords with respect to the $s$-blocks. Given
$p\in\{0\ld a_{s+1}-1\}$, we let
\[
\cA_p \ = \ \{m\in\Z\mid m=\nu a_{s+1} + 
  p \textrm{ for some } \nu\in\Z\} \ 
\]
and write $p(m)=p$ if $m\in\cA_p$.  The second parameter $q$ describes the
position with respect to the nearest $n$-block, where $n$ is chosen such that
$a_n/2 < N \leq a_{n+1}/2$, and the third parameter $\iota$ determines the local
configuration of symbols around the $n$-block. Both are somewhat more subtle to
define. 

Suppose first that the interval $I^{m,N}=\{m-N\ld m+N\}$ does not
intersect any $n$-block. Then we let $q(m)=0$ and
$\iota(m)=0$. Otherwise, $I^{m,N}$ intersects exactly one $n$-block
$B(m)=\{-a_n\ld a_n\}+\ell a_{n+1}$. In this case, we define $j(m) = \ell
a_{n+1}-m$ and let
\[
q(m) \ = q \quad \textrm{if} \quad j(m) \in
\left(\left.\frac{qN}{2^{j_0}},\frac{(q+1)N}{2^{j_0}}\right.\right]
\]
Note that $j(m)$ lies between $-N-a_n$ and $N+a_n$ and thus $q(m)$
ranges only from at least $-3\cdot 2^{j_0}$ to at most $3\cdot
2^{j_0}-1$.

In order to define $\iota$ for the case $I^{m,N}$ intersects an $n$-block
$B(m)$, we denote by $B^-(m)$ and $B^+(m)$ the two neareast $n$-blocks to the
left, respectively, right of $B(m)$. Further, we denote by $J^-(m)$ the interval
between $B^-(m)$ and $B(m)$ and by $J^+(m)$ the interval between $B(m)$ and
$B^+(m)$. Note that by construction, there exist unique integers $n^\pm$ such
that $J^\pm(m)\smin A_n\ssq C_{n^\pm}$. Hence, $\omega_k$ remains constant on
each of the sets $J^\pm(m)\smin A_n$. We assume that $n$ is even, such that
$\omega_k=0$ for all $k\in B(m)\smin A_{n-1}=B(m)\cap C_n$. Then there are four
possibilities:
\begin{itemize}
\item[(1)] $\omega_k=0$ for all $k\in (J^-(m)\cup J^+(m))\smin A_n$;
\item[(2)] $\omega_k=0$ for all $k\in J^-(m)\smin A_n$ and $\omega_k=1$
  for all $k\in J^+(m)\smin A_n$;
\item[(3)] $\omega_k=1$ for all $k\in J^-(m)\smin A_n$ and $\omega_k=0$
  for all $k\in J^+(m)\smin A_n$;
\item[(4)] $\omega_k=1$ for all $k\in (J^-(m)\cup J^+(m))\smin A_n$.
\end{itemize}
We let $\iota(m)=\iota$ if case $(\iota)$ applies. Now, given $q\in\{-3\cdot
2^{j_0}\ld 3\cdot 2^{j_0}-1\}$ and $\iota\in\{0\ld 4\}$ we let
\[
\cA^\iota_{p,q} \ = \ \left\{m\in\cA_p \ \mid \ q(m)=q \textrm{ and }
  \iota(m)=\iota \right\} \ ,
\]
where we set $\cA^0_{p,q}=\emptyset$ if $q\neq 0$. This defines the required
decomposition of $\Z$ into at most $30\cdot 2^{j_0}\cdot a_{s+1}$ partition
elements. It remains to show that given $m,m'\in \cA^\iota_{p,q}$, the words
$\omega^{m,N}$ and $\omega^{m',N}$ cannot be $\delta$-separated with respect to
$D_N$.  Thus, we need to estimate the maximal number of mismatches that can
appear between two such words.

First, since the position of the words with respect to the
$a_{s+1}$-periodic set $A_s$ is identical, we have that $j+m\in A_s$
if and only if $j+m'\in A_s$, and in this case $\omega^{m,N}_j=
\omega^{m',N}_j$. 
 
Secondly, if either $j+m\in A_{n-1}\smin A_s$ or $j+m'\in A_{n-1}\smin A_s$,
then this might result in a mismatch. However, since
\[
\frac{\#\left( I^{m,N} \cap A_{n-1}\right)}{\# I^{m,N}} \ \leq \
2\sum_{i=s+1}^{n-1} b^{-1}_i \ < \ \delta/4 \ 
\]
and likewise for $I^{m',N}$, we have that the contribution of such
mismatches to $D_N(\omega^{m,N},\omega^{m',N})$ is at most $\delta/2$.

Finally, it remains to count the possible mismatches with $j+m$, $j+m'\notin
A_{n-1}$. If $j(m)=j(m')$, then there are no such mismatches, since the
intervals $B(m)-m$ and $B(m')-m'$ as well as $J^\pm(m)-m$ and $J^\pm(m')-m'$
coincide and case $(\iota)$ above applies to both $m$ and $m'$. Otherwise, there
are possible overlaps between non-corresponding intervals, but since
$|j(m)-j(m')|\leq N/2^{j_0}$ these overlaps concern at most $2N/2^{j_0} <
2N\cdot (\delta/4)$ positions. Again, this results in a contribution to
$D_N(\omega^{m,N},\omega^{m',N})$ of at most $\delta/2$. Altogether, this
yields $D_N(\omega^{m,N},\omega^{m',N})<\delta$ as required and thus completes
the proof.\qed\medskip

\begin{remark}\mbox{}\alphlist
\item
  Using similar, but simpler arguments, it is possible to show that the power
  entropy of the above example equals 1.
\item At the same time, it can be shown that the amorphic complexity
  of the example is infinite, and even $S^*_\nu(f,\delta)=\infty$ for
  sufficiently small $\nu,\delta>0$. This is, in fact, a consequence
  of a much more general statement. It is possible to prove that the
  asymptotic separation numbers $S^*_\nu(f,\delta)$ of a minimal
  action of a homeomorphism $f$ on a compact metric space are all finite if
  and only if the system is Weyl mean equicontinuous. By a recent result of
  Downarowicz and Glasner, this holds if and only if the system is an
  isomorphic extension of its maximal equicontinuous factor \cite[Theorem
  2.1]{DownarowiczGlasner2015}. In particular, it needs to be uniquely
  ergodic, which is not true for our example due to the fluctuating
  symbol frequencies. Since these issues will be explored further in
  \cite{GroegerJaeger2014AmorphicComplexityQuasicrystals}, we do not
  go into further detail here.
\item A construction which is very similar to the above one, but results in a
  uniquely ergodic irregular Toeplitz flow, can be found in
  \cite{downarowicz2015odometers}. Analogous arguments can be applied to show
  that this example also has modified power entropy zero.\listend  
\end{remark}

\section{Conclusions and open questions}\label{Conclusions}

It may seem, admittedly, that this note has a somewhat negative touch, since the
presented results are mostly negative ones. We mainly showed that modified power
entropy does not satisfy a variational principle, is not independent of
transient dynamics, does not respond to the transition from equicontinuous
systems to their almost 1-1 extensions and cannot be used either to distinguish
between regular and irregular extensions of minimal equicontinuous
systems. However, as we discussed in Section~\ref{Discussion} already, these
issues do not have to be seen as disadvantages of the notion itself. As said
before, the existence of a variational principle and the insensitivity to
transient effects are not necessarily positive features of a slow entropy, since
this depends very much on the purpose one has in mind.  Thus, the presented
facts should rather be understood as clarifications and simply imply that for
the specific aspects we concentrate on other topological invariants have to be
identified in order to fulfill the respective tasks or requirements. We also
note that it follows from results of Ferenczi that modified power entropy does
detect the transition from uniquely ergodic isomorphic to non-isomorphic
extensions of compact group rotations
\cite{Ferenczi1997MeasureTheoreticComplexity}. 

As we have seen, the transition from equicontinuous minimal systems to
their almost 1-1 extensions can be detected by means of $a$-entropy or
amorphic complexity (with suitably chosen scale functions). In the
other cases, however, this leads to the following open questions.
\begin{itemize}
\item[(a)] Is there a topological invariant $h$ for continuous maps on (compact)
  metric or topological spaces that gives meaningful information about zero
  entropy systems, but at the same time satisfies $h(f)=h(f_{|\Omega(f)})$?
\item[(b)] Is there such a topological invariant that satisfies a variational
  principle with respect to a suitable measure-theoretic analogue?
\end{itemize}
We note that Kong and Chen \cite{KongChen2014} recently introduced a slow entropy
which satisfies a `non-standard' variational principle, in which the supremum is
taken over all probability measures on the phase space (and not just the
invariant ones).
\begin{itemize}\item[(c)] Is there a meaningful topological invariant for dynamical systems
  which is zero for all regular almost 1-1 extensions of
  equicontinuous systems, but strictly positive for all irregular
  almost 1-1 extensions of such systems?
\end{itemize}
Some progress on closely related questions has recently been made by Li, Tu and
Ye \cite{LiTuYe2014} and Downarowicz and Glasner \cite{DownarowiczGlasner2015}.

\end{document}